\documentclass[12pt]{article}
\usepackage[francais,english]{babel}

\usepackage[]{amsfonts}
\usepackage{amssymb}
\usepackage{amsmath}

\def\couleur(#1 #2 #3)
	{
\special{ps::[end]}}

\def\underset#1#2{\mathrel{\mathop{\kern0pt #2}\limits_{#1}}}

\def\overset#1#2{\mathrel{\mathop{\kern0pt #2}\limits^{#1}}}

\def\bx#1{\setbox1=\hbox{\kern3pt{#1}\kern3pt}			
 \dimen1=\ht1 \advance\dimen1 by 3pt \dimen2=\dp1 \advance\dimen2 by 3pt
 \setbox1=\hbox{\vrule height\dimen1 depth\dimen2\box1\vrule}%
 \setbox1=\vbox{\hrule\box1\hrule}%
 \advance\dimen1 by .4pt \ht1=\dimen1
 \advance\dimen2 by .4pt \dp1=\dimen2 \box1\relax}

\def\wbb#1{\kern#1em}
\def\vci{\vrule  width.02em height1.47ex depth-.0ex}		
\def\11{{\rm\wbb{.2}\vci\wbb{-.37}1}}

\newcommand{\Supp}{\operatorname{Supp}} 

\newtheorem{Theorem}{Theorem}[section]
\newtheorem{Lemma}[Theorem]{Lemma}
\newtheorem{Remark}[Theorem]{Remark}

\begin{document}

\title{The raising steps method. Applications to the  $\bar \partial $  equation in Stein manifolds.}

\author{Eric Amar}

\date{}
\maketitle
 \ \par 
\renewcommand{\abstractname}{Abstract}

\begin{abstract}
In order to get estimates on the solutions of the equation 
 $\bar \partial u=\omega $  on Stein manifold, we introduce a
 new method the "raising steps method", to get global results
 from local ones. In particular it allows us to transfer results
 form open sets in  ${\mathbb{C}}^{n}$  to open sets in a Stein manifold.\ \par 
Using it we get  $\displaystyle L^{r}-L^{s}$  results for solutions
 of equation  $\bar \partial u=\omega $  with a gain,  $\displaystyle
 s>r,$  in strictly pseudo convex domains in Stein manifolds.\ \par 
We also get  $\displaystyle L^{r}-L^{s}$  results for domains
 in  ${\mathbb{C}}^{n}$  locally biholomorphic to convex domains
 of finite type.\ \par 
\end{abstract}

\section{Introduction.}
\quad  In all the sequel a domain of a smooth manifold  $X$  will be
 a connected open set  $\displaystyle \Omega $  of  $X,$  relatively
 compact and with  ${\mathcal{C}}^{\infty }$ smooth boundary.\ \par 
\quad \quad  	A strictly pseudo convex domain  $\displaystyle \Omega $  of
 the complex manifold  $X,$   $\displaystyle s.p.c.$  for short,
 is a domain defined by a smooth real valued function  $\rho
 $  on  $X$  such that\ \par 
\quad \quad 	 $\bullet $   $\rho $  is strictly pseudo convex near  $\displaystyle
 \bar \Omega \ ;$ \ \par 
\quad \quad 	 $\bullet $   $\displaystyle \Omega :=\lbrace z\in X::\rho (z)<0\rbrace
 \ ;$ \ \par 
\quad \quad 	 $\bullet $   $\displaystyle \forall z\in \partial \Omega ,\
 \partial \rho (z)\neq 0.$ \ \par 
\quad \quad  	The starting point of this work is a result of S. Krantz~\cite{KrantzDbar76}
 for  $\displaystyle (0,1)$  forms and generalized to  $\displaystyle
 (0,q)$  forms in~\cite{AmarLrLsCnv13} :\ \par 
\begin{Theorem}
 ~\label{rS13}Let  $\displaystyle \Omega $  be a  $\displaystyle
 s.p.c.$  domain in  ${\mathbb{C}}^{n}$  then for  $\displaystyle
 1<r<2n+2$  we have\par 
\quad \quad 	 $\bullet $  	 	 $\displaystyle \forall \omega \in L^{r}_{(p,q)}(\Omega
 ),\ \bar \partial \omega =0,\ \exists u\in L^{s}_{(p,q-1)}(\Omega
 )::\bar \partial u=\omega ,\ {\left\Vert{u}\right\Vert}_{L^{s}(\Omega
 )}\lesssim {\left\Vert{\omega }\right\Vert}_{L^{r}(\Omega )},$ \par 
with  $\displaystyle \ \frac{1}{s}=\frac{1}{r}-\frac{1}{2(n+1)}.$ \par 
\quad  $\bullet $ 	 For  $\displaystyle r=2n+2$  we have\par 
\quad \quad \quad \quad \quad   $\displaystyle \exists u\in BMO_{(p,q)}(\Omega )::\bar \partial
 u=\omega ,\ {\left\Vert{u}\right\Vert}_{BMO(\Omega )}\lesssim
 {\left\Vert{\omega }\right\Vert}_{L^{2n+2}(\Omega )}.$ \par 
\quad  If  $\omega $  is a  $\displaystyle (p,1)$  form we have also :\par 
\quad  $\bullet $ 	 for  $\displaystyle r=1,$ \par 
\quad \quad \quad \quad \quad   $\displaystyle \exists u\in L^{s,\infty }_{(p,0)}(\Omega )::\bar
 \partial u=\omega ,\ {\left\Vert{u}\right\Vert}_{L^{s,\infty
 }(\Omega )}\lesssim {\left\Vert{\omega }\right\Vert}_{L^{1}(\Omega )}$ \par 
with  $\displaystyle \ \frac{1}{s}=1-\frac{1}{2(n+1)}.$ \par 
\quad  $\bullet $ 	 for  $\displaystyle r>2n+2,$ \par 
\quad \quad \quad \quad \quad   $\displaystyle \exists u\in \Gamma ^{\beta }_{(p,0)}(\Omega
 )::\bar \partial u=\omega ,\ {\left\Vert{u}\right\Vert}_{\Gamma
 ^{\beta }(\Omega )}\lesssim {\left\Vert{\omega }\right\Vert}_{L^{r}(\Omega
 )},$ \par 
where  $\displaystyle \beta =1-\frac{2n+2}{r}$  and  $\displaystyle
 \ \Gamma ^{\beta }$  is an anisotropic Lipschitz class of functions.\par 
Moreover the solution  $u$  is linear on the data  $\omega .$ \par 
\end{Theorem}
\quad \quad  	These results strongly improve the  $\displaystyle L^{r}-L^{r}$
  results of N. Ovrelid~\cite{Ovrelid71} because they give a
 gain,  $\displaystyle s>r,$  and this is the key for the "raising
 steps method" to work. \ \par 
\quad  This method allows 	to get from local results on solutions 
 $u$  of an equation  $\displaystyle Du=\omega ,$  global ones
 in a smooth manifold  $X.$ \ \par 
\quad \quad  	The applications will be done in the several complex variables
 setting so 	on any complex manifold  $X$  we define first the
 "Lebesgue measure" as in H\"ormander's book~\cite{Hormander73}
 section 5.2, with a hermitian metric locally equivalent to the
 usual one on any analytic coordinates patch. Associated to this
 metric there is a volume form  $dm$  and we take it for the
 Lebesgue measure on  $X.$ \ \par 
\quad \quad  	There already exist results on solutions of the  $\bar \partial
 $  equation in  $\displaystyle s.p.c.$  domains in complex manifolds.
 For instance H\"ormander~\cite{Hormander73} solve the  $\bar
 \partial $  equation with  $\displaystyle L^{2}-L^{2}$  estimates
 for any  $\displaystyle (p,q)$  currents. Kerzman~\cite{Kerzman71}
 proved a  $\displaystyle L^{r}-L^{r}$  estimate on solution
 of  $\bar \partial $  for any  $\displaystyle r\in \lbrack 1,\infty
 \lbrack $  and a H\"older estimate in case  $\displaystyle r=\infty
 $  for  $\displaystyle (0,1)$  currents. In~\cite{HenkinLeiterer84},
 chap. 4, Henkin and Leiterer built global kernels on  $\displaystyle
 s.p.c.$  domains in a Stein manifold for  $\displaystyle (0,q)$
  forms and get precise uniform estimates. Demailly and Laurent~\cite{DemaillyLaurent87},
 built global kernels on  $\displaystyle s.p.c.$  domains in
 a Stein manifold for  $\displaystyle (p,q)$  forms and get 
  $\displaystyle L^{r}-L^{r}$  estimate on solution of  $\bar
 \partial $  for any  $\displaystyle r\in \lbrack 1,\infty \lbrack
 $  in the case the Stein manifold is equipped with a hermitian
 metric with null curvature.\ \par 
\ \par 
\quad \quad  	As a first application of this "raising steps" method we get
 the following theorems which improve all already known results.\ \par 
\begin{Theorem}
 ~\label{lrEstimees34}Let  $\displaystyle \Omega $  be a smoothly
 bounded relatively compact strictly pseudo convex domain in
 the Stein manifold  $\displaystyle X.$  \par 
 There is a constant  $\displaystyle C=C(p,q,r)>0$  such that
 if  $\omega \in L^{r}_{(p,q)}(\Omega ),\ \bar \partial \omega
 =0$  if  $\displaystyle 1\leq r\leq 2,$  there is a  $\displaystyle
 (p,q-1)$  current  $u$  such that, with  $\displaystyle \gamma
 :=\min (\frac{1}{2(n+1)},\frac{1}{r}-\frac{1}{2}),$  and  $\displaystyle
 \ \frac{1}{s}=\frac{1}{r}-\gamma ,$ \par 
\quad \quad \quad \quad \quad 	 $\bar \partial u=\omega $  and  $\displaystyle u\in L^{s}_{(p,q-1)}(\Omega
 ),\ {\left\Vert{u}\right\Vert}_{L_{(p,q-1)}^{s}(\Omega )}\leq
 C{\left\Vert{\omega }\right\Vert}_{L_{(p,q)}^{r}(\Omega )}.$ \par 
\end{Theorem}
\quad \quad  	Let  ${\mathcal{H}}_{p}^{r}(\Omega )$  be the set of all  $(p,\
 0)\ \bar \partial $  closed forms in  $\Omega $  and in  $\displaystyle
 L^{r}(\Omega ).$  In the sequel we shall denote  $\displaystyle
 r'$  the conjugate exponent of  $\displaystyle r,\ \frac{1}{r}+\frac{1}{r'}=1.$
 \ \par 
To deal with  $\displaystyle r>2$  we have to make the assumption
 that  $\displaystyle \omega \perp {\mathcal{H}}_{p}^{r'}(\Omega
 )$  if  $\displaystyle q=n$  and that  $\omega $  has a {\sl
 compact support} for  $\displaystyle q<n,$   i.e.  $\displaystyle
 \omega \in L^{r,c}_{(p,q)}(\Omega )$  if  $\displaystyle q<n.$ \ \par 
\begin{Theorem}
 ~\label{lrEstimees33}Let  $\displaystyle \Omega $  be a smoothly
 bounded relatively compact spc domain in the Stein manifold
  $\displaystyle X.$  Let  $\displaystyle \ \frac{2(n+1)}{n+2}\leq
 r<2(n+1),$  there is a constant  $\displaystyle C=C_{r}>0$ 
 such that if   $\omega \in L^{r,c}_{(p,q)}(\Omega ),\ \bar \partial
 \omega =0$  if  $\displaystyle 1\leq q<n$  and  $\displaystyle
 \omega \in L^{r}_{(p,q)}(\Omega ),\ \omega \perp {\mathcal{H}}_{p}^{r'}(\Omega
 )$  if  $\displaystyle q=n,$  then there is a  $\displaystyle
 (p,q-1)$  current  $u$  such that, with  $\displaystyle \ \frac{1}{s}=\frac{1}{r}-\frac{1}{2(n+1)},$
 \par 
\quad \quad \quad \quad \quad 	 $\bar \partial u=\omega $  and  $\displaystyle u\in L^{s}_{(p,q-1)}(\Omega
 ),\ {\left\Vert{u}\right\Vert}_{s}\leq C{\left\Vert{\omega }\right\Vert}_{r}.$
 \par 
\end{Theorem}
\quad \quad  	In the case of  $\displaystyle (0,1)$  currents, these assumptions
 are not necessary, thanks to the results of N. Kerzman~\cite{Kerzman71}.\ \par 
\begin{Theorem}
 Let  $\displaystyle \Omega $  be a smoothly bounded relatively
 compact spc domain in the Stein manifold  $\displaystyle X.$ \par 
\quad  $\bullet $ 	 If  $\displaystyle 1\leq r<2(n+1)$  there is a
 constant  $\displaystyle C=C_{r}>0$  such that  if   $\omega
 \in L^{r}_{(0,1)}(\Omega ),\ \bar \partial \omega =0$  then
 there is a function  $u$  such that, with  $\displaystyle \
 \frac{1}{s}=\frac{1}{r}-\frac{1}{2(n+1)},$  \par 
\quad \quad \quad \quad \quad 	 $\bar \partial u=\omega $  and  $\displaystyle u\in L^{s}(\Omega
 ),\ {\left\Vert{u}\right\Vert}_{L^{s}(\Omega )}\leq C{\left\Vert{\omega
 }\right\Vert}_{L_{(0,1)}^{r}(\Omega )}.$ \par 
\quad \quad 	 $\bullet $  If  $\displaystyle r=2(n+1)$  then for any  $\displaystyle
 s<\infty $  there is a function  $u$  such that\par 
\quad \quad \quad \quad \quad \quad 	 $\bar \partial u=\omega $  and  $\displaystyle u\in L^{s}(\Omega
 ),\ {\left\Vert{u}\right\Vert}_{L^{s}(\Omega )}\leq C{\left\Vert{\omega
 }\right\Vert}_{L_{(0,1)}^{r}(\Omega )}.$ \par 
\quad \quad 	 $\bullet $  if  $\displaystyle 2(n+1)<r<\infty $  there is
 a function  $u$  such that\par 
\quad \quad \quad \quad \quad \quad \quad   $\displaystyle u\in \Lambda ^{\beta }(\Omega )::\bar \partial
 u=\omega ,\ {\left\Vert{u}\right\Vert}_{\Lambda ^{\beta }(\Omega
 )}\leq C{\left\Vert{\omega }\right\Vert}_{L_{(0,1)}^{r}(\Omega )},$ \par 
where  $\displaystyle \beta =\frac{1}{2}-\frac{(n+1)}{r}$  and
  $\displaystyle \ \Lambda ^{\beta }$  is the H\"older class
 of functions of order  $\beta .$ \par 
\end{Theorem}
\quad \quad  	Another way to release the compact support assumption is to
 have a Stein manifold equipped with a hermitian  metric of null
 curvature in the sense of Demailly Laurent~\cite{DemaillyLaurent87},
 because then we can use their  $\displaystyle L^{r}-L^{r}$ 
 estimates for any  $\displaystyle r\in \lbrack 1,\infty \lbrack .$ \ \par 
\begin{Theorem}
 Let  $\displaystyle \Omega $  be a smoothly bounded relatively
 compact spc domain in the Stein manifold  $\displaystyle X$
  equipped with a metric with null curvature. \par 
\quad  $\bullet $ 	 If  $\displaystyle 1\leq r<2(n+1)$  and  $\omega
 \in L^{r}_{(p,q)}(\Omega ),\ \bar \partial \omega =0.$  Then
 there is a  $\displaystyle (p,q-1)$  current  $u$  such that\par 
\quad \quad \quad \quad \quad 	 $\bar \partial u=\omega $  and  $\displaystyle u\in L^{s}_{(p,q-1)}(\Omega
 ),\ {\left\Vert{u}\right\Vert}_{s}\leq C{\left\Vert{\omega }\right\Vert}_{L_{(p,q)}^{r}(\Omega
 )}$  with  $\displaystyle \ \frac{1}{s}=\frac{1}{r}-\frac{1}{2(n+1)}.$ \par 
\quad \quad 	 $\bullet $  If  $\displaystyle r=2(n+1)$   and  $\omega \in
 L^{r}_{(p,q)}(\Omega ),\ \bar \partial \omega =0.$  Then for
 any  $\displaystyle s<\infty $  there is a  $\displaystyle (p,q-1)$
  current  $u$  such that\par 
\quad \quad \quad \quad \quad 	 $\bar \partial u=\omega $  and  $\displaystyle u\in L^{s}_{(p,q-1)}(\Omega
 ),\ {\left\Vert{u}\right\Vert}_{L^{s}_{(p,q-1)}}\leq C{\left\Vert{\omega
 }\right\Vert}_{L_{(p,q)}^{r}(\Omega )}.$ \par 
\end{Theorem}
\quad \quad  	We notice that N. Kerzman in~\cite{Kerzman71}, in order to
 solve  $\bar \partial $  for  $\displaystyle (0,1)$  forms,
 also use local solutions and to find a global one he used also
 the H\"ormander  $\displaystyle L^{2}$  solution, but his method
 is based on "bump" around point at the boundary and is completely
 different from the raising steps method introduced here.\ \par 
\quad \quad  	As another application of this method we get the following theorem :\ \par 
\begin{Theorem}
 Let  $\displaystyle \Omega $  be a domain in  ${\mathbb{C}}^{n}$
  such that locally around any point of  $\displaystyle \partial
 \Omega ,$   $\displaystyle \Omega $  is biholomorphic to a convex
 domain of finite type at most of type  $m\ ;$  then if  $\omega
 \in L^{r}_{(p,q)}(\Omega ),\ \bar \partial \omega =0$  if  $\displaystyle
 1\leq r\leq 2,$  there is a  $\displaystyle (p,q-1)$  current
  $u$  such that, with  $\displaystyle \gamma :=\min (\frac{1}{mn+2},\frac{1}{r}-\frac{1}{2}),$
  and  $\displaystyle \ \frac{1}{s}=\frac{1}{r}-\gamma ,$ \par 
\quad \quad \quad \quad \quad 	 $\bar \partial u=\omega $  and  $\displaystyle u\in L^{s}_{(p,q-1)}(\Omega
 ),\ {\left\Vert{u}\right\Vert}_{L_{(p,q-1)}^{s}(\Omega )}\leq
 C{\left\Vert{\omega }\right\Vert}_{L_{(p,q)}^{r}(\Omega )}.$ \par 
\end{Theorem}
\quad \quad  	To deal with  $\displaystyle r>2$  we have again to make the
 assumption that  $\displaystyle \omega \perp {\mathcal{H}}_{p}^{r'}(\Omega
 )$  if  $\displaystyle q=n$  and that  $\omega $  has a compact
 support for  $\displaystyle q<n,$   i.e.  $\displaystyle \omega
 \in L^{r,c}_{(p,q)}(\Omega )$  if  $\displaystyle q<n.$ \ \par 
\begin{Theorem}
Let  $\displaystyle \Omega $  be a domain in  ${\mathbb{C}}^{n}$
  such that locally around any point of  $\displaystyle \partial
 \Omega ,$   $\displaystyle \Omega $  is biholomorphic to a convex
 domain of finite type at most of type  $m\ ;$  Let  $\displaystyle
 \ 2<r\leq mn+2$  and  $\omega \in L^{r,c}_{(p,q)}(\Omega ),\
 \bar \partial \omega =0$  if  $\displaystyle 1\leq q<n$  and
  $\displaystyle \omega \in L^{r}_{(p,q)}(\Omega ),\ \omega \perp
 {\mathcal{H}}_{p}^{r'}(\Omega )$  if  $\displaystyle q=n\ ;$
  then there is a  $\displaystyle (p,q-1)$  current  $u$  such
 that, with  $\displaystyle \ \frac{1}{s}=\frac{1}{r}-\frac{1}{mn+2},$ \par 
\quad \quad \quad \quad \quad 	 $\bar \partial u=\omega $  and  $\displaystyle u\in L^{s}_{(p,q-1)}(\Omega
 ),\ {\left\Vert{u}\right\Vert}_{s}\leq C{\left\Vert{\omega }\right\Vert}_{r}.$
 \par 
\end{Theorem}

\section{The "raising steps" method.}
\quad \quad  	We shall deal with the following situation : we have a  ${\mathcal{C}}^{\infty
 }$ smooth manifold  $X$  admitting partitions of unity and a
 decreasing scale  $\displaystyle \lbrace B_{r}\rbrace _{r\geq
 1},\ s\geq r\Rightarrow \ B_{s}\subset B_{r}$  of Banach spaces
 of functions or forms defined on open sets of  $X.$  These Banach
 spaces must be modules over  ${\mathcal{D}},$  the space of
  $\displaystyle {\mathcal{C}}^{\infty }$  functions with compact
 support, i.e.\ \par 
\quad \quad 	 $\displaystyle \forall \Omega \ open\ in\ X,\ \forall \chi
 \in {\mathcal{D}}(\Omega ),\ \exists C(\chi )>0::\forall f\in
 B_{r},\ \chi f\in B_{r}(\Omega )$  and  $\displaystyle \ {\left\Vert{\chi
 f}\right\Vert}_{B_{r}}\leq C(\chi ){\left\Vert{f}\right\Vert}_{B_{r}}.$ \ \par 
For instance  $\displaystyle B_{r}=L^{r}$  the Lebesgue spaces,
 or  $\displaystyle B_{r}=H_{r}^{2}$  the Sobolev spaces, etc...\ \par 
\quad  We are interested in solution of the linear equation  $\displaystyle
 Du=\omega ,$  where  $D$  is a linear operator, with eventually
 the constraint  $\displaystyle \Delta \omega =0,$  where  $\Delta
 $  is also a linear operator such that  $\displaystyle \Delta
 ^{2}=0.$  In case there is no constraint we take  $\displaystyle
 \Delta \equiv 0.$  One aim is to apply this to the  $\bar \partial
 $  equation.\ \par 
\quad \quad  	We shall put the following hypotheses on  $D,$  for any domain
  $\displaystyle \Omega \subset X$  :\ \par 
\quad \quad  	(i)  $\displaystyle \forall \chi \in {\mathcal{D}}(\Omega ),\
 D\chi \in {\mathcal{D}}(\Omega )\ ;$ \ \par 
\quad \quad  	(ii)   $\displaystyle \forall \chi \in {\mathcal{D}}(\Omega
 ),\ \forall f\in B_{r}(\Omega ),\ D(\chi f)=D\chi \cdot f+\chi Df\ ;$ \ \par 
as can be easily seen, a linear differential operator  $D$ 
 verifies these assumptions.\ \par 
\quad \quad  	Let  $\displaystyle \Omega $  be a relatively compact domain
 in  $X.$  Now we shall make the following assumptions on  $X$
  and  $\displaystyle \Omega .$  There is a  $\displaystyle r_{0}>1$
  and a   $\displaystyle \delta >0$  such that, setting  $\displaystyle
 \ \frac{1}{s}=\frac{1}{r}-\delta ,$ \ \par 
\quad \quad  	(iii) there is a covering  $\displaystyle \lbrace U_{j}\rbrace
 _{j=1,...,N}$  of  $\displaystyle \bar \Omega $  such that,
  $\displaystyle \forall r\leq r_{0},$  if  $\displaystyle \omega
 \in B_{r}(\Omega ),\ \Delta \omega =0,$  we can solve  $Du_{j}=\omega
 $  in  $\displaystyle U_{j}\cap \Omega $  with  $\displaystyle
 B_{r}(\Omega )-B_{s}(\Omega \cap U_{j})$  estimates, i.e.  $\displaystyle
 \exists C_{0}>0$  such that\ \par 
\quad \quad \quad \quad \quad 	 $\displaystyle \exists u_{j}\in B_{s}(\Omega _{j}),\ Du_{j}=\omega
 $  in  $\displaystyle \Omega _{j}:=U_{j}\cap \Omega $  and 
 $\displaystyle \ {\left\Vert{u_{j}}\right\Vert}_{B_{s}(\Omega
 _{j})}\leq C_{0}{\left\Vert{\omega }\right\Vert}_{B_{r}(\Omega )}.$ \ \par 
\quad \quad  	(iv) We can solve  $Dw=\omega $  globally in  $\displaystyle
 \Omega $  with  $\displaystyle B_{r_{0}}-B_{r_{0}}$  estimates, i.e.\ \par 
\quad \quad \quad \quad \quad 	 $\displaystyle \exists E>0,\ \exists w::Dw=\omega $  in  $\displaystyle
 \Omega $  and  $\displaystyle \ {\left\Vert{w}\right\Vert}_{B_{r_{0}}(\Omega
 )}\leq E{\left\Vert{\omega }\right\Vert}_{B_{r_{0}}(\Omega )}$
  provided that  $\displaystyle \Delta \omega =0.$ \ \par 
\quad \quad  	Then we have\ \par 
\begin{Theorem}
 ~\label{LrStrPsStein37}(Raising steps theorem) Under the assumptions
 above, there is a constant  $\displaystyle C>0,$  for  $\displaystyle
 r\leq r_{0},$  if  $\omega \in B_{r}(\Omega ),\ \Delta \omega
 =0$  there is a $\displaystyle \ u\in B_{s}(\Omega )$  with
  $\displaystyle \gamma :=\min (\delta ,\frac{1}{r}-\frac{1}{r_{0}}),$
  and  $\displaystyle \ \frac{1}{s}=\frac{1}{r}-\gamma ,$  such that\par 
\quad \quad \quad \quad \quad 	 $Du=\omega $  and  $\displaystyle u\in B_{s}(\Omega ),\ {\left\Vert{u}\right\Vert}_{B_{s}(\Omega
 )}\leq C{\left\Vert{\omega }\right\Vert}_{B_{r}(\Omega )}.$ \par 
\end{Theorem}
\quad \quad  	Proof.\ \par 
Let  $\displaystyle r\leq r_{0}$  and  $\displaystyle \omega
 \in B_{r}(\Omega ),\ \Delta \omega =0$  ; we start with the
 covering  $\displaystyle \lbrace U_{j}\rbrace _{j=1,...,N}$
  and the local solution  $Du_{j}=\omega $  with  $\displaystyle
 u_{j}\in B_{s}(\Omega _{j})$  given by hypothesis (iii).\ \par 
Let  $\displaystyle \chi _{j}$  be a  ${\mathcal{C}}^{\infty
 }$ smooth partition of unity subordinate to  $\displaystyle
 \lbrace U_{j}\rbrace _{j=1,...,N}$  and set\ \par 
\quad \quad \quad \quad \quad 	 $\displaystyle v_{0}:=\sum_{j=1}^{N}{\chi _{j}u_{j}}.$ \ \par 
Then we have\ \par 
\quad \quad 	 $\bullet $   $\displaystyle v_{0}\in B_{s_{0}}(\Omega )$  because
  $\displaystyle B_{s}$  is a module over  $\displaystyle {\mathcal{D}}(\Omega
 )$  and  $\displaystyle \ {\left\Vert{v_{0}}\right\Vert}_{B_{s_{0}}(\Omega
 )}\leq C{\left\Vert{\omega }\right\Vert}_{B_{r}(\Omega )}$ 
 with  $\displaystyle \ \frac{1}{s_{0}}=\frac{1}{r}-\delta $
  and  $\displaystyle C=NC_{0}\max _{j=1,...,N}C(\chi _{j})$
  by hypothesis (iii).\ \par 
\quad \quad 	 $\bullet $   $\displaystyle Dv_{0}=\sum_{j=1}^{N}{\chi _{j}Du_{j}}+\sum_{j=1}^{N}{D\chi
 _{j}\wedge u_{j}}$  by hypothesis (ii) hence\ \par 
\quad \quad \quad \quad \quad   $\displaystyle Dv_{0}=\sum_{j=1}^{N}{\chi _{j}\omega }+\sum_{j=1}^{N}{D\chi
 _{j}\wedge u_{j}}=\omega +\omega _{1}$ \ \par 
with\ \par 
\quad \quad \quad \quad \quad 	 $\displaystyle \omega _{1}:=\sum_{j=1}^{N}{D\chi _{j}\wedge u_{j}}.$ \ \par 
But, by hypothesis (i),  $\displaystyle D\chi _{j}\in {\mathcal{D}}(\Omega
 )$  hence because  $\displaystyle B_{s}$  is a module over 
 $\displaystyle {\mathcal{D}}(\Omega ),$  we have  $\displaystyle
 \omega _{1}\in B_{s_{0}}(\Omega ),$  with  $\displaystyle \
 \frac{1}{s_{0}}=\frac{1}{r}-\delta $  and  $\displaystyle \
 {\left\Vert{\omega _{1}}\right\Vert}_{B_{s_{0}}(\Omega )}\leq
 G{\left\Vert{\omega }\right\Vert}_{B_{r}(\Omega )}$  with  $G=C_{0}N\max
 _{j=1,...,N}C(D\chi _{j}).$ \ \par 
Hence the regularity of  $\omega _{1}$  is higher of one step
  $\delta $  than that of  $\omega .$ \ \par 
Moreover  $\Delta \omega _{1}=\Delta ^{2}v_{0}-\Delta \omega
 =0$  because  $\Delta ^{2}=0$  and  $\Delta \omega =0.$ \ \par 
\quad \quad  	Two cases\ \par 
Case 1 : if  $\displaystyle \ \frac{1}{s_{0}}=\frac{1}{r}-\delta
 \Rightarrow s_{0}\geq r_{0}$  we can solve a global  $D$  in
  $\displaystyle B_{r_{0}}(\Omega ),$  by assumption (iv) i.e.\ \par 
\quad \quad \quad \quad \quad 	 $\displaystyle \exists w\in B_{r_{0}}(\Omega )::Dw=\omega _{1},\
 {\left\Vert{w}\right\Vert}_{B_{r_{0}}(\Omega )}\leq E{\left\Vert{\omega
 _{1}}\right\Vert}_{B_{r_{0}}(\Omega )}\Rightarrow {\left\Vert{w}\right\Vert}_{B_{r_{0}}(\Omega
 )}\leq EG{\left\Vert{\omega }\right\Vert}_{B_{r}(\Omega )}.$ \ \par 
It remains to set\ \par 
\quad \quad \quad \quad \quad 	 $\displaystyle u:=v_{0}-w$ \ \par 
to have, by the linearity of  $D,$ \ \par 
\quad \quad \quad \quad \quad 	 $Du=\omega $ \ \par 
and, with  $\displaystyle \ \frac{1}{s}=\frac{1}{r}-\gamma ,\
 \gamma :=\min (\delta ,\frac{1}{r}-\frac{1}{r_{0}}),$  and 
 $\displaystyle F:=EG,$ \ \par 
\quad \quad \quad \quad \quad 	 $\displaystyle u\in B_{s}(\Omega ),\ {\left\Vert{u}\right\Vert}_{B_{s}(\Omega
 )}\leq F{\left\Vert{\omega }\right\Vert}_{B_{r}(\Omega )}.$ \ \par 
\ \par 
Case 2 : if  $\displaystyle \ \frac{1}{s_{0}}=\frac{1}{r}-\delta
 \Rightarrow s_{0}<r_{0}$   then we continue :\ \par 
\quad \quad \quad \quad \quad 	 $\exists v_{1}\in B_{s_{1}}(\Omega ),\ Dv_{1}=\omega _{1}+\omega _{2}$ \ \par 
with\ \par 
\quad \quad \quad \quad \quad 	 $\displaystyle \ {\left\Vert{v_{1}}\right\Vert}_{B_{s_{1}}(\Omega
 )}\leq C{\left\Vert{\omega _{1}}\right\Vert}_{B_{s_{0}}(\Omega
 )}\leq CG{\left\Vert{\omega }\right\Vert}_{B_{r}(\Omega )},$ \ \par 
\quad \quad \quad \quad \quad   $\displaystyle \ \omega _{2}\in B_{s_{1}}(\Omega ),\ {\left\Vert{\omega
 _{2}}\right\Vert}_{B_{s_{1}}(\Omega )}\leq G{\left\Vert{\omega
 _{1}}\right\Vert}_{B_{s_{0}}(\Omega )},$ \ \par 
\quad \quad \quad \quad \quad 	 $\displaystyle \ \frac{1}{s_{1}}=\frac{1}{s_{0}}-\delta =\frac{1}{r}-2\delta
 ,$ \ \par 
\quad \quad \quad \quad \quad 	 $\Delta \omega _{2}=0.$ \ \par 
And\ \par 
\quad \quad \quad \quad \quad 	 $D(v_{0}-v_{1})=\omega +\omega _{1}-Dv_{1}=\omega -\omega _{2}.$ \ \par 
Hence by induction :\ \par 
\quad \quad \quad \quad \quad 	 $\displaystyle \exists v_{0},...,v_{N}::v_{j}\in B_{s_{j}}(\Omega
 ),\ \frac{1}{s_{j}}=\frac{1}{r}-(j+1)\delta ,\ {\left\Vert{v_{j}}\right\Vert}_{B_{s_{j}}(\Omega
 )}\leq CG^{j-1}{\left\Vert{\omega }\right\Vert}_{B_{r}(\Omega )}$ \ \par 
and\ \par 
\quad \quad \quad \quad \quad 	 $\displaystyle D(\sum_{j=0}^{N}{(-1)^{j}v_{j}})=\omega +(-1)^{N}\omega
 _{N}$ \ \par 
with\ \par 
\quad \quad \quad \quad \quad   $\displaystyle \omega _{N}\in B_{s_{N-1}}(\Omega ),\ {\left\Vert{\omega
 _{N}}\right\Vert}_{B_{s_{N-1}}(\Omega )}\leq G^{N}{\left\Vert{\omega
 }\right\Vert}_{B_{r}(\Omega )},\ \frac{1}{s_{N-1}}=\frac{1}{r}-N\delta
 $ \ \par 
and  $\displaystyle \Delta \omega _{N}=0.$ \ \par 
\quad \quad  	Hence the regularity of  $\omega _{N}$  raises of  $N$  steps
  $\delta $  from that of  $\omega .$ \ \par 
\quad  Now let  $\displaystyle s_{0}::\frac{1}{s_{0}}=\frac{1}{r}-\delta
 $  and suppose that  $\displaystyle s_{0}\leq r_{0}$  then take
  $N$  such that  $\displaystyle \ \frac{1}{s_{N-1}}=\frac{1}{r}-N\delta
 \leq \frac{1}{r_{0}}$  then  $\displaystyle \omega _{N}\in B_{s_{N-1}}(\Omega
 )\subset B_{r_{0}}(\Omega )$  and now we can solve a {\sl global}
  $D$  in  $\displaystyle B_{r_{0}}(\Omega ),$  by assumption (iv) i.e.\ \par 
\quad \quad \quad \quad \quad 	 $\displaystyle \exists w\in B_{r_{0}}(\Omega )::Dw=\omega _{N},\
 {\left\Vert{w}\right\Vert}_{B_{r_{0}}(\Omega )}\leq E{\left\Vert{\omega
 _{N}}\right\Vert}_{B_{r_{0}}(\Omega )}\Rightarrow {\left\Vert{w}\right\Vert}_{B_{r_{0}}(\Omega
 )}\leq ECG^{N}{\left\Vert{\omega }\right\Vert}_{B_{r}(\Omega )}.$ \ \par 
It remains to set\ \par 
\quad \quad \quad \quad \quad 	 $\displaystyle u:=\sum_{j=0}^{N}{(-1)^{j}v_{j}}+(-1)^{N}w$ \ \par 
to have\ \par 
\quad \quad \quad \quad \quad 	 $Du=\omega $ \ \par 
and, with  $\displaystyle \ \frac{1}{s}=\frac{1}{r}-\delta $
  and  $\displaystyle F:=ECG^{N},$ \ \par 
\quad \quad \quad \quad \quad 	 $\displaystyle u\in B_{s}(\Omega ),\ {\left\Vert{u}\right\Vert}_{B_{s}(\Omega
 )}\leq F{\left\Vert{\omega }\right\Vert}_{B_{r}(\Omega )}.$
   $\blacksquare $ \ \par 

\section{Applications.}
\quad \quad  	We proved in~\cite{AmarLrLsCnv13} the following theorem for
 strictly pseudo convex (s.p.c.) domains in  ${\mathbb{C}}^{n}.$ \ \par 
\begin{Theorem}
 ~\label{raisingSteps11}Let  $\displaystyle \Omega $  be a  $\displaystyle
 s.p.c.$  domain in  ${\mathbb{C}}^{n}$  then\par 
\quad \quad 	 $\bullet $  for  $\displaystyle 1<r<2n+2$  we have\par 
\quad \quad \quad \quad \quad 	 $\displaystyle \forall \omega \in L^{r}_{(p,q)}(\Omega ),\
 \bar \partial \omega =0,\ \exists u\in L^{s}_{(p,q-1)}(\Omega
 )::\bar \partial u=\omega ,\ {\left\Vert{u}\right\Vert}_{L^{s}(\Omega
 )}\lesssim {\left\Vert{\omega }\right\Vert}_{L_{(p,q)}^{r}(\Omega )},$ \par 
with  $\displaystyle \ \frac{1}{s}=\frac{1}{r}-\frac{1}{2(n+1)}.$ \par 
\quad  $\bullet $ 	 For  $\displaystyle r=2n+2$  we have\par 
\quad \quad \quad \quad \quad   $\displaystyle \exists u\in \bigcap_{r\geq 1}{L^{r}_{(p,q-1)}(\Omega
 )}::\bar \partial u=\omega .$ \par 
If  $\omega $  is a  $\displaystyle (p,1)$  form we have also :\par 
\quad  $\bullet $ 	 for  $\displaystyle r=1,$  we have\par 
\quad \quad \quad \quad \quad   $\displaystyle \exists u\in L^{s,\infty }_{(p,0)}(\Omega )::\bar
 \partial u=\omega ,\ {\left\Vert{u}\right\Vert}_{L^{s,\infty
 }(\Omega )}\lesssim {\left\Vert{\omega }\right\Vert}_{L^{1}(\Omega )}$ \par 
with  $\displaystyle \ \frac{1}{s}=1-\frac{1}{2(n+1)}.$ \par 
\quad  $\bullet $ 	 for  $\displaystyle r>2n+2$  we have\par 
\quad \quad \quad \quad \quad   $\displaystyle \exists u\in \Lambda ^{\beta }(\Omega )::\bar
 \partial u=\omega ,\ {\left\Vert{u}\right\Vert}_{\Lambda ^{\beta
 }(\Omega )}\lesssim {\left\Vert{\omega }\right\Vert}_{L^{r}(\Omega )},$ \par 
where  $\displaystyle \beta =\frac{1}{2}-\frac{(n+1)}{r}$  and
  $\displaystyle \ \Lambda ^{\beta }$  is the H\"older class
 of functions of order  $\beta .$ \par 
Moreover the solution  $u$  is linear on the data  $\omega .$ \par 
\end{Theorem}
\quad \quad  	In fact the theorem we have gave results on  $\displaystyle
 BMO(\Omega )$  instead of  $\displaystyle \ \bigcap_{r\geq 1}{L^{r}_{(p,q-1)}(\Omega
 )}$  with control of the norm and on an anisotropic H\"older
 class  $\displaystyle \Gamma ^{\beta }(\Omega )$  instead of
 the usual H\"older class. But because we are mainly interested
 in classical results, it is not necessary to go further here.
 The aim is to have this kind of results with a Stein manifold
 instead of  ${\mathbb{C}}^{n}.$ \ \par 

\subsection{The transfer from  ${\mathbb{C}}^{n}$  to Stein manifold.}
\quad \quad  	We shall apply the raising steps method to the case of  $D=\bar
 \partial ,\ \Delta =\bar \partial ,\ B_{r}=L_{(p,q)}^{r}$  the
 space of  $\displaystyle (p,q)$  currents with coefficients
 in the Lebesgue space  $\displaystyle L^{r}$  and  $\displaystyle
 X$  a Stein manifold. Clearly (i) and (ii) are verified.\ \par 
\begin{Lemma}
~\label{LrStrPsStein49}	Let  $\displaystyle \Omega $  be a 
 $\displaystyle s.p.c.$  domain in the complex manifold  $\displaystyle
 X.$  There is a covering  $\displaystyle \lbrace (U_{j},\varphi
 _{j})\rbrace _{j=1,...,N}$  of  $\displaystyle \bar \Omega $
  by coordinates patches of  $X$  such that  $\displaystyle \Omega
 _{j}:=U_{j}\cap \Omega $  is still relatively compact,  ${\mathcal{C}}^{\infty
 }$  smoothly bounded and strictly pseudo convex in  $\displaystyle X.$ \par 
\end{Lemma}
\quad \quad  	Proof.\ \par 
Let  $\displaystyle z\in \bar \Omega $  and  $\displaystyle
 (G_{z},\varphi _{z})$  a coordinates patch such that  $\displaystyle
 z\in G_{z}.$  Set  $V_{z}:=\varphi _{z}(G_{z})\subset {\mathbb{C}}^{n}$
  and  $\zeta :=\varphi _{z}(z).$ \ \par 
We have two cases \ \par 
\quad \quad  	(i) if  $\displaystyle z\in \Omega $  take a ball  $B_{\zeta
 }$  centered at  $\zeta $  and contained in  $\displaystyle
 V_{z}.$   $B$  is spc in  ${\mathbb{C}}^{n}\ ;$  set  $\displaystyle
 U_{z}:=\varphi _{z}^{-1}(B)$  then  $\displaystyle U_{z}=U_{z}\cap
 \Omega $  is still spc with smooth boundary in  $X,$  because
  $\varphi _{z}$  is biholomorphic in  $\displaystyle G_{z},$
  and  $\displaystyle U_{z}$  is a neighbourhood of  $\displaystyle z.$ \ \par 
\quad \quad  	(ii) if  $\displaystyle z\in \partial \Omega $  we look at
  $V_{z}\subset {\mathbb{C}}^{n}$  and we make the completion
 of the part of the boundary of  $\displaystyle \varphi _{z}(G_{z}\cap
 \Omega )$  in  $\displaystyle V_{z}$  to get a smoothly bounded
 strictly pseudo convex domain  $\Gamma _{z}$  in  $V_{z}$  as
 in~\cite{AmCoh84}.\ \par 
\quad \quad  	Now we extend  $\displaystyle \Gamma _{z}$  on the other side
 of  $\displaystyle \varphi _{z}(G_{z}\cap \partial \Omega )$
  in  $\displaystyle V_{z}$  as an open set  $W_{z}$  such that
  $\displaystyle W_{z}\cap \varphi _{z}(U_{z}\cap \Omega )=\Gamma
 _{z}.$  Now  we set  $U_{z}:=\varphi _{z}^{-1}(W_{z}),$  then
  $\displaystyle \Omega _{z}:=U_{z}\cap \Omega =\varphi _{z}^{-1}(\Gamma
 _{z})$  is strictly pseudo convex with smooth boundary in  $X.$ \ \par 
\quad \quad  	Hence we have that  $\displaystyle \lbrace U_{z}\rbrace _{z\in
 \bar \Omega }$  is a covering of the compact set  $\displaystyle
 \bar \Omega $  and we can extract of it a finite subset  $\displaystyle
 \lbrace U_{j}\rbrace _{j=1,...,N}$  which is still a covering
 of  $\displaystyle \bar \Omega $  with all the required properties.
  $\blacksquare $ \ \par 
\ \par 
\begin{Remark}
 ~\label{raisingSteps212}We have that the Lebesgue measure on
  $X$  restricted to  $\displaystyle U_{j}$  is equivalent to
 the restriction of the Lebesgue measure of  ${\mathbb{C}}^{n}$
  to  $\displaystyle \varphi _{j}(U_{j}).$  This equivalence
 is uniform with respect to  $\displaystyle j=1,...,N$  of course.
 Moreover we have that the distance in  $\displaystyle \Omega
 \cap U_{j}$  is also uniformly equivalent to the distance in
  $\Gamma _{j}$  because  $\varphi _{j}$  is biholomorphic in
  $\displaystyle U_{j}$  and there is only a finite number of
  $\displaystyle U_{j}$  to deal with.\par 
\end{Remark}

\subsection{Case with use of the  $\displaystyle L^{2}$  estimates
 of H\"ormander.}
\quad \quad  	Let  $\displaystyle \Omega $  be  $\displaystyle s.p.c.$  domain
 in the Stein manifold  $\displaystyle X.$  We plan to apply
 the raising steps method to extend estimates from domains in
  ${\mathbb{C}}^{n}$  to domains in  $\displaystyle X.$ \ \par 
\quad \quad  	By lemma~\ref{LrStrPsStein49} we have a covering  $\displaystyle
 \lbrace U_{j}\rbrace _{j=1,...,N}$  of  $\displaystyle \bar
 \Omega $  by coordinates patches  $\displaystyle (U_{j},\varphi
 _{j})$  such that\ \par 
\quad \quad \quad \quad \quad 	 $\displaystyle \Omega _{j}=U_{j}\cap \Omega $  is spc in  $X$
  and  $\displaystyle \Gamma _{j}=\varphi _{j}(\Omega _{j})$
  is spc in  ${\mathbb{C}}^{n}$  with  ${\mathcal{C}}^{\infty
 }$  smooth boundary.\ \par 
Let  $\displaystyle \omega \in L_{(p,q)}^{r}(\Omega ),\ $ we
 have by remark~\ref{raisingSteps212}, that Lebesgue measures
 on  $\displaystyle \Omega _{j}$  and on  $\displaystyle \Gamma
 _{j}$  are equivalent,\ \par 
\quad \quad \quad \quad \quad 	 $\displaystyle \omega \in L_{(p,q)}^{r}(\Omega )$  implies
  $\displaystyle \varphi _{j}^{*}\omega \in L_{(p,q)}^{r}(\Gamma _{j}).$ \ \par 
\quad \quad  	Hence we can apply theorem~\ref{raisingSteps11} to each  $\displaystyle
 \Gamma _{j}$  to get a  $u'_{j}\in L^{s}_{(p,q-1)}(\Gamma _{j}),\
 \bar \partial u'_{j}=\varphi _{j}^{*}\omega $  with  $\displaystyle
 \ \frac{1}{s}=\frac{1}{r}-\frac{1}{2(n+1)}.$  So we have here
  $\displaystyle \delta =\frac{1}{2(n+1)}.$ \ \par 
Back to  $X,$  we have a  $\displaystyle u_{j}\in L^{s}_{(p,q-1)}(\Omega
 _{j}),\ \bar \partial u_{j}=\omega $  in  $\displaystyle U_{j}\cap
 \Omega =\Omega _{j}$  with control of the norm. So assumption
 (iii) is fulfilled.\ \par 
\quad  We get the assumption (iv) by the well known theorem of H\"ormander~\cite{Hormander73},
 ~\cite{FolKoh72}.\ \par 
\begin{Theorem}
 ~\label{lrEstimees35}Let  $\displaystyle \Omega $  be a  $\displaystyle
 s.p.c.$  domain in the Stein manifold  $\displaystyle X.$  There
 is a constant  $\displaystyle C>0$  such that if  $\omega \in
 L^{2}_{(p,q)}(\Omega ),\ \bar \partial \omega =0$  there is
 a  $\displaystyle (p,q-1)$  current  $u\in L^{2}_{(p,q-1)}(\Omega
 ),\ \bar \partial u=\omega $  and  $\ {\left\Vert{u}\right\Vert}_{2}\leq
 C{\left\Vert{\omega }\right\Vert}_{2}.$ \par 
\end{Theorem}
\quad \quad  	So an application of the "raising steps" theorem~\ref{LrStrPsStein37}
 gives\ \par 
\begin{Theorem}
 ~\label{LrStrPsStein38}Let  $\displaystyle \Omega $  be a  $\displaystyle
 s.p.c.$  domain in the Stein manifold  $\displaystyle X.$  There
 is a constant  $\displaystyle C>0$  such that if  $\omega \in
 L^{r}_{(p,q)}(\Omega ),\ \bar \partial \omega =0,$  with  $\displaystyle
 r\leq 2,$ \par 
set  $\displaystyle \gamma :=\min (\frac{1}{2(n+1)},\frac{1}{r}-\frac{1}{2})$
  and  $\displaystyle \ \frac{1}{s}=\frac{1}{r}-\gamma ,$  then
 there is a  $\displaystyle (p,q-1)$  current  $u$  such that\par 
\quad \quad \quad \quad \quad 	 $\bar \partial u=\omega $  and  $\displaystyle u\in L^{s}_{(p,q-1)}(\Omega
 ),\ {\left\Vert{u}\right\Vert}_{L_{(p,q-1)}^{s}(\Omega )}\leq
 C{\left\Vert{\omega }\right\Vert}_{L_{(p,q)}^{r}(\Omega )}.$ \par 
\end{Theorem}
\ \par 
\quad \quad  	To deal with the case  $\displaystyle r>2$  we shall proceed
 by duality and ask that  $\omega $  has {\sl compact support},
 i.e.  $\displaystyle \omega \in L^{r,c}_{(p,q)}(\Omega )$  
 when  $\displaystyle q<n.$ \ \par 
Recall that  ${\mathcal{H}}_{p}^{r}(\Omega )$  is the set of
 all  $(p,\ 0)\ \bar \partial $  closed forms in  $\Omega $ 
 and in  $\displaystyle L^{r}(\Omega ).$  As usual let  $\displaystyle
 r'$  the conjugate exponent to  $\displaystyle r,\ \frac{1}{r}+\frac{1}{r'}=1.$
 \ \par 
\begin{Theorem}
 ~\label{lrEstimees33}Let  $\displaystyle \Omega $  be a  $\displaystyle
 s.p.c.$  domain in the Stein manifold  $\displaystyle X.$  Let
  $\displaystyle \ \frac{2(n+1)}{n+2}\leq r\leq 2(n+1)$  and
  $\omega \in L^{r,c}_{(p,q)}(\Omega ),\ \bar \partial \omega
 =0$  if  $\displaystyle 1\leq q<n$  and  $\displaystyle \omega
 \in L^{r}_{(p,q)}(\Omega ),\ \omega \perp {\mathcal{H}}_{p}^{r'}(\Omega
 )$  if  $\displaystyle q=n\ ;$  then there is a  $\displaystyle
 (p,q-1)$  current  $u$  such that, with  $\displaystyle \ \frac{1}{s}=\frac{1}{r}-\frac{1}{2(n+1)},$
 \par 
\quad \quad \quad \quad \quad 	 $\bar \partial u=\omega $  and  $\displaystyle u\in L^{s}_{(p,q-1)}(\Omega
 ),\ {\left\Vert{u}\right\Vert}_{s}\leq C{\left\Vert{\omega }\right\Vert}_{r}.$
 \par 
\end{Theorem}
\quad \quad  	Proof.\ \par 
We use the same technique of duality we already use in~\cite{AnGrauLr12}
 inspired by the Serre duality theorem~\cite{Serre55}.\ \par 
\begin{Lemma}
 ~\label{AndreGrau189} For  $\displaystyle \Omega ,\ \omega $
  as in the theorem, consider the function  ${\mathcal{L}}={\mathcal{L}}_{\omega
 },$  defined on  $(n-p,n-q+1)$  form  $\alpha \in L^{s'}(\Omega
 ),\ \bar \partial $  closed in  $\Omega ,$  as follows:\par 
\quad \quad \quad \quad \quad 	 ${\mathcal{L}}(\alpha ):=(-1)^{p+q-1}{\left\langle{\omega ,\varphi
 }\right\rangle},$  where  $\varphi \in L^{r'}(\Omega )$  is
 such that  $\bar \partial \varphi =\alpha $  in  $\displaystyle \Omega .$ \par 
Then  ${\mathcal{L}}$ {\it  }is well defined and linear.\par 
\end{Lemma}
\quad \quad  	Proof.\ \par 
First notice that if  $\displaystyle \ \frac{1}{s}=\frac{1}{r}-\frac{1}{2(n+1)}$
  then  $\displaystyle \ \frac{1}{r'}=\frac{1}{s'}-\frac{1}{2(n+1)}.$
  Such a current  $\varphi $  with  $\bar \partial \varphi =\alpha
 $  exists since  $\displaystyle s',$  the conjugate exponent
 of  $s$  verifies  $\displaystyle s'\leq 2,$  hence we can apply
 theorem~\ref{LrStrPsStein38}, with the remark that this time
  $\displaystyle s'\leq r'.$ \ \par 
Suppose first that  $q<n.$ \ \par 
In order for  ${\mathcal{L}}$  to be well defined we need\ \par 
\quad \quad \quad   $\forall \varphi ,\psi \in L^{r'}_{(n-p,n-q)}(\Omega ),\ \bar
 \partial \varphi =\bar \partial \psi \Rightarrow {\left\langle{\omega
 ,\varphi }\right\rangle}={\left\langle{\omega ,\psi }\right\rangle}.$ \ \par 
This is meaningful because  $\omega \in L^{r,c}(\Omega ),\ r>1,\
 \Supp \omega \Subset \Omega .$ \ \par 
Then we have  $\bar \partial (\varphi -\psi )=0$  hence  $\displaystyle
 \varphi -\psi \in L^{r'}(\Omega )\subset L^{s'}(\Omega )$  so
 we can solve  $\bar \partial $  in  $L^{r'}(\Omega )$  because
  $\displaystyle s'\leq 2$  by theorem~\ref{LrStrPsStein38} :\ \par 
\quad \quad \quad \quad \quad 	 $\displaystyle \exists \gamma \in L^{r'}_{(n-p,n-q-1)}(\Omega
 )::\bar \partial \gamma =(\varphi -\psi )\in L^{r'}(\Omega ).$ \ \par 
So  $\ {\left\langle{\omega ,\varphi -\psi }\right\rangle}={\left\langle{\omega
 ,\bar \partial \gamma }\right\rangle}=(-1)^{p+q-1}{\left\langle{\bar
 \partial \omega ,\gamma }\right\rangle}=0$  because  $\omega
 $  being compactly supported in  $\Omega $  there is no boundary term.\ \par 
Hence  ${\mathcal{L}}$  is well defined in that case.\ \par 
\quad \quad  	Suppose now that  $q=n.$ \ \par 
Of course  $\bar \partial \omega =0$  but we have that  $\varphi
 ,\ \psi $  are  $(p,\ 0)$  forms hence  $\bar \partial (\varphi
 -\psi )=0$  means that  $h:=\varphi -\psi $  is a  $\bar \partial
 $  closed  $(p,\ 0)$  form hence  $h\in {\mathcal{H}}_{p}(\Omega
 ).$  The hypothesis  $\displaystyle \omega \perp {\mathcal{H}}_{p}^{r'}(\Omega
 )$  gives  $\ {\left\langle{\omega ,h}\right\rangle}=0,$  and
  ${\mathcal{L}}$  is also well defined in that case. We notice
 that  $\omega $  with compact support is not needed in that
 case, but in order that the scalar product be defined, we need
  $\displaystyle h\in {\mathcal{H}}_{p}^{r'}(\Omega ).$ \ \par 
\quad \quad  	It remains to see that  ${\mathcal{L}}$  is linear, so let
  $\alpha =\alpha _{1}+\alpha _{2},$  with  $\displaystyle \alpha
 _{j}\in L^{s'}(\Omega ),\ \bar \partial \alpha _{j}=0,\ j=1,2\
 ;$  we have  $\alpha =\bar \partial \varphi ,\ \alpha _{1}=\bar
 \partial \varphi _{1}$  and  $\alpha _{2}=\bar \partial \varphi
 _{2},$  with  $\varphi ,\ \varphi _{1},\ \varphi _{2}$  in 
 $\displaystyle L^{r'}(\Omega )$  so, because  $\bar \partial
 (\varphi -\varphi _{1}-\varphi _{2})=0,$  we have\ \par 
\quad \quad  	if  $\displaystyle q<n\ :$ \ \par 
\quad \quad \quad \quad \quad 	 $\varphi =\varphi _{1}+\varphi _{2}+\bar \partial \psi ,$ 
 with  $\psi $  in  $L^{r'}(\Omega ),$  as we did above,\ \par 
so\ \par 
\quad \quad \quad \quad \quad 	 ${\mathcal{L}}(\alpha )=(-1)^{p+q-1}{\left\langle{\omega ,\varphi
 }\right\rangle}=(-1)^{p+q-1}{\left\langle{\omega ,\varphi _{1}+\varphi
 _{2}+\bar \partial \psi }\right\rangle}=\ {\mathcal{L}}(\alpha
 _{1})+{\mathcal{L}}(\alpha _{2})+(-1)^{p+q-1}{\left\langle{\omega
 ,\bar \partial \psi }\right\rangle},$ \ \par 
but  $\ {\left\langle{\omega ,\bar \partial \psi }\right\rangle}={\left\langle{\bar
 \partial \omega ,\psi }\right\rangle}=0,$  because  $\displaystyle
 \Supp \omega \Subset \Omega $  implies there is no boundary term.\ \par 
Hence  ${\mathcal{L}}(\alpha )={\mathcal{L}}(\alpha _{1})+{\mathcal{L}}(\alpha
 _{2}).$ \ \par 
\quad \quad  	if  $\displaystyle q=n\ :$ \ \par 
	because  $\bar \partial (\varphi -\varphi _{1}-\varphi _{2})=0,$
  we have  $\displaystyle h:=\varphi -\varphi _{1}-\varphi _{2}\in
 {\mathcal{H}}_{p}^{r'}(\Omega )$  and the hypothesis  $\displaystyle
 \omega \perp {\mathcal{H}}_{p}^{r'}(\Omega )$  gives  $\ {\left\langle{\omega
 ,h}\right\rangle}=0,$  so\ \par 
\quad \quad \quad \quad \quad 	 $\displaystyle {\mathcal{L}}(\alpha )=(-1)^{p+q-1}{\left\langle{\omega
 ,\varphi }\right\rangle}=(-1)^{p+q-1}{\left\langle{\omega ,\varphi
 _{1}+\varphi _{2}+h}\right\rangle}=\ {\mathcal{L}}(\alpha _{1})+{\mathcal{L}}(\alpha
 _{2})+(-1)^{p+q-1}{\left\langle{\omega ,h}\right\rangle},$ \ \par 
hence  $\displaystyle {\mathcal{L}}(\alpha )={\mathcal{L}}(\alpha
 _{1})+{\mathcal{L}}(\alpha _{2}).$ \ \par 
We notice again that  $\omega $  with compact support is not
 needed in that case.\ \par 
\quad \quad  	The same for  $\alpha =\lambda \alpha _{1}$  and the linearity.
  $\blacksquare $ \ \par 
\begin{Lemma}
~\label{lrEstimees32}Still with the same hypotheses as above
 there is a  $\displaystyle (p,q-1)$  current  $u$  such that\par 
\quad \quad 	 $\displaystyle \forall \alpha \in L_{(n-p,n-q+1)}^{s'}(\Omega
 ),\ {\left\langle{u,\alpha }\right\rangle}={\mathcal{L}}(\alpha
 )=(-1)^{p+q-1}{\left\langle{\omega ,\varphi }\right\rangle},$ \par 
and\par 
\quad \quad \quad \quad \quad 	 $\displaystyle \sup _{\alpha \in L^{s'}(\Omega ),\ {\left\Vert{\alpha
 }\right\Vert}_{L^{s'}(\Omega )}\leq 1}\ \left\vert{{\left\langle{u,\alpha
 }\right\rangle}}\right\vert \leq C{\left\Vert{\omega }\right\Vert}_{L^{r}(\Omega
 )}.$ \par 
\end{Lemma}
\quad \quad  	Proof.\ \par 
By lemma~\ref{AndreGrau189} we have that  ${\mathcal{L}}$  is
 a linear form on  $(n-p,n-q+1)$  forms  $\alpha \in L^{s'}(\Omega
 ),\ \bar \partial $  closed in  $\Omega .$ \ \par 
	We have\ \par 
\quad \quad \quad \quad \quad 	 $\exists \varphi \in L_{(n-p,n-q)}^{r'}(\Omega )::\bar \partial
 \varphi =\alpha ,\ {\left\Vert{\varphi }\right\Vert}_{L^{r'}(\Omega
 )}\leq C_{2}{\left\Vert{\alpha }\right\Vert}_{L^{s'}(\Omega )}$ \ \par 
and by its very definition\ \par 
\quad \quad \quad \quad \quad 	 ${\mathcal{L}}(\alpha )={\left\langle{\omega ,\ \varphi }\right\rangle}.$
 \ \par 
	By H\"older inequalities\ \par 
\quad \quad \quad \quad \quad 	 $\ \left\vert{{\mathcal{L}}(\alpha )}\right\vert \leq {\left\Vert{\omega
 }\right\Vert}_{L^{r}(\Omega )}{\left\Vert{\varphi }\right\Vert}_{L^{r'}(\Omega
 )}={\left\Vert{\omega }\right\Vert}_{L^{r}(\Omega )}{\left\Vert{\varphi
 }\right\Vert}_{L^{r'}(\Omega )}.$ \ \par 
But there is a constant  $C$  such that\ \par 
\quad \quad \quad \quad \quad 	 $\ {\left\Vert{\varphi }\right\Vert}_{L^{r'}(\Omega )}\leq
 C{\left\Vert{\alpha }\right\Vert}_{L^{s'}(\Omega )}.$ \ \par 
Hence\ \par 
\quad \quad \quad \quad \quad 	 $\ \left\vert{{\mathcal{L}}(\alpha )}\right\vert \leq C{\left\Vert{\omega
 }\right\Vert}_{L^{r}(\Omega )}{\left\Vert{\alpha }\right\Vert}_{L^{s'}(\Omega
 )}.$ \ \par 
\quad \quad  	So we have that the norm of  ${\mathcal{L}}$  is bounded on
 the subspace of  $\bar \partial $  closed forms in  $L^{s'}(\Omega
 )$  by  $\ C{\left\Vert{\omega }\right\Vert}_{L^{r}(\Omega )}.$ \ \par 
\quad \quad  	We apply the Hahn-Banach theorem to extend  ${\mathcal{L}}$
  with the {\sl same} norm to {\sl all}   $(n-p,n-q+1)$  forms
 in  $L^{s'}(\Omega ).$  As in Serre duality theorem (~\cite{Serre55},
 p. 20) this is one of the main ingredient in the proof.\ \par 
\quad \quad  	This means, by the definition of currents, that there is a
  $(p,q-1)$  current  $u$  which represents the extended form
  ${\mathcal{L}},$  i.e.\ \par 
\quad \quad \quad \quad \quad 	 $\displaystyle \forall \alpha \in L_{(n-p,n-q+1)}^{s'}(\Omega
 ),\ {\left\langle{u,\alpha }\right\rangle}={\mathcal{L}}(\alpha
 )=(-1)^{p+q-1}{\left\langle{\omega ,\varphi }\right\rangle},$ \ \par 
and such that\ \par 
\quad \quad \quad \quad \quad 	 $\displaystyle \sup _{\alpha \in L^{s'}(\Omega ),\ {\left\Vert{\alpha
 }\right\Vert}_{L^{s'}(\Omega )}\leq 1}\ \left\vert{{\left\langle{u,\alpha
 }\right\rangle}}\right\vert \leq C{\left\Vert{\omega }\right\Vert}_{L^{r}(\Omega
 )}.$ \ \par 
 $\blacksquare $ \ \par 
\quad \quad  	The lemma~\ref{lrEstimees32} says that\ \par 
\quad \quad \quad \quad \quad 	 $\displaystyle \forall \alpha \in L_{(n-p,n-q+1)}^{s'}(\Omega
 ),\ {\left\langle{u,\alpha }\right\rangle}={\mathcal{L}}(\alpha
 )=(-1)^{p+q-1}{\left\langle{\omega ,\varphi }\right\rangle},$ \ \par 
hence applied to  $\displaystyle \varphi \in {\mathcal{D}}_{(n-p,n-q)}(\Omega
 )$  we get  $\alpha =\bar \partial \varphi \in {\mathcal{D}}_{(n-p,n-q+1)}\subset
 L^{s'}(\Omega )$  and\ \par 
\quad \quad \quad \quad \quad 	 $\ {\left\langle{u,\bar \partial \varphi }\right\rangle}={\mathcal{L}}(\bar
 \partial \varphi )=(-1)^{p+q-1}{\left\langle{\omega ,\varphi
 }\right\rangle}\Rightarrow {\left\langle{\bar \partial u,\varphi
 }\right\rangle}={\left\langle{\omega ,\varphi }\right\rangle}$ \ \par 
because  $\varphi $  has compact support in  $\displaystyle
 \Omega ,$  hence  $\bar \partial u=\omega $  in the sense of
 distributions.\ \par 
Moreover the fact that\ \par 
\quad \quad \quad \quad \quad 	 $\displaystyle \sup _{\alpha \in L^{s'}(\Omega ),\ {\left\Vert{\alpha
 }\right\Vert}_{L^{s'}(\Omega )}\leq 1}\ \left\vert{{\left\langle{u,\alpha
 }\right\rangle}}\right\vert \leq C{\left\Vert{\omega }\right\Vert}_{L^{r}(\Omega
 )}$ \ \par 
gives, by an easy duality, that\ \par 
\quad \quad \quad \quad \quad 	 $\displaystyle \ {\left\Vert{u}\right\Vert}_{L_{(p,q-1)}^{s}(\Omega
 )}\leq C{\left\Vert{\omega }\right\Vert}_{L_{(p,q)}^{r}(\Omega )}$ \ \par 
which completes the proof of theorem~\ref{lrEstimees33}.  $\blacksquare
 $ \ \par 
\begin{Remark}
 The theorem is valid up to  $\displaystyle r=2(n+1)$  because
 there we have to take  $\displaystyle s'=1$  and we get\par 
\quad \quad \quad \quad \quad 	 $\displaystyle \sup _{\alpha \in L^{1}(\Omega ),\ {\left\Vert{\alpha
 }\right\Vert}_{L^{1}(\Omega )}\leq 1}\ \left\vert{{\left\langle{u,\alpha
 }\right\rangle}}\right\vert \leq C{\left\Vert{\omega }\right\Vert}_{L^{r}(\Omega
 )}$ \par 
which implies that  $\displaystyle u\in L_{(p,q-1)}^{\infty
 }(\Omega )$  with control of the norm.\par 
\end{Remark}

\subsection{Case of  $\displaystyle (0,1)$  forms by use of
 Kerzman's estimates.}
\quad \quad  	In the special case of  $\displaystyle (0,1)$  forms  $\omega
 ,$  we apply the raising steps method with Kerzman's estimates~\cite{Kerzman71}
 for the global solution.\ \par 
\quad \quad  	Let  $\displaystyle \Omega $  be a  $\displaystyle s.p.c.$
  domain in the Stein manifold  $\displaystyle X.$  As above
 we have the assumption (iii) fulfilled with  $\displaystyle
 \delta =\frac{1}{2(n+1)}.$  We shall prove\ \par 
\begin{Theorem}
 Let  $\displaystyle \Omega $  be a  $\displaystyle s.p.c.$ 
 domain in the Stein manifold  $\displaystyle X.$  \par 
\quad  $\bullet $ 	 If  $\displaystyle 1\leq r<2(n+1)$  there is a
 constant  $\displaystyle C=C_{r}>0$  such that  if   $\omega
 \in L^{r}_{(0,1)}(\Omega ),\ \bar \partial \omega =0$  then
 there is a function  $u$  such that, with  $\displaystyle \
 \frac{1}{s}=\frac{1}{r}-\frac{1}{2(n+1)},$  \par 
\quad \quad \quad \quad \quad 	 $\bar \partial u=\omega $  and  $\displaystyle u\in L^{s}(\Omega
 ),\ {\left\Vert{u}\right\Vert}_{L^{s}(\Omega )}\leq C{\left\Vert{\omega
 }\right\Vert}_{L_{(0,1)}^{r}(\Omega )}.$ \par 
\quad \quad 	 $\bullet $  If  $\displaystyle r=2(n+1)$  then for any  $\displaystyle
 s<\infty $  there is a function  $u$  such that\par 
\quad \quad \quad \quad \quad \quad 	 $\bar \partial u=\omega $  and  $\displaystyle u\in L^{s}(\Omega
 ),\ {\left\Vert{u}\right\Vert}_{L^{s}(\Omega )}\leq C{\left\Vert{\omega
 }\right\Vert}_{L_{(0,1)}^{r}(\Omega )}.$ \par 
\quad \quad 	 $\bullet $  if  $\displaystyle r>2(n+1)$  there is a function
  $u$  such that\par 
\quad \quad \quad \quad \quad \quad \quad   $\displaystyle u\in \Lambda ^{\beta }_{(p,0)}(\Omega )::\bar
 \partial u=\omega ,\ {\left\Vert{u}\right\Vert}_{\Lambda ^{\beta
 }(\Omega )}\lesssim {\left\Vert{\omega }\right\Vert}_{L_{(0,1)}^{r}(\Omega
 )},$ \par 
where  $\displaystyle \beta =\frac{1}{2}-\frac{(n+1)}{r}$  and
  $\displaystyle \ \Lambda ^{\beta }$  is the H\"older class
 of functions of order  $\beta .$ \par 
\end{Theorem}
\quad \quad  	Proof.\ \par 
The assumption (iv) is true for any  $\displaystyle r_{0}\in
 \lbrack 1,\infty \rbrack $  by Kerzman's estimates~\cite{Kerzman71}
 so starting with a  $\displaystyle r\in \lbrack 1,2(n+1)\lbrack
 ,$  we choose  $\displaystyle r_{0}=s$  with  $\displaystyle
 \ \frac{1}{s}=\frac{1}{r}-\frac{1}{2(n+1)}>0,$  hence an application
 of the raising steps theorem~\ref{LrStrPsStein37} gives the first point.\ \par 
\quad  For the second point the proof of the raising steps theorem
 gives that we have  $\displaystyle v_{0}:=\sum_{j=1}^{N}{\chi
 _{j}u_{j}}$  with  $\displaystyle u_{j}\in \bigcap_{1\leq t}{L^{t}(\Omega
 \cap U_{j})}$  hence  $\displaystyle v_{0}\in \bigcap_{1\leq
 t}{L^{t}(\Omega \cap U_{j})}$  because the  $\chi _{j}$  are
 in  $\displaystyle {\mathcal{D}}(\Omega \cap U_{j}).$  Choose
 a  $s<\infty $  and  $\displaystyle r_{0}=s,$  then apply Kerzman's
 result : there is a correction  $\displaystyle w\in L^{s}(\Omega
 )$  such that  $\bar \partial (v_{0}-w)=\omega $  and  $\displaystyle
 v_{0}-w\in L^{s}(\Omega ).$ \ \par 
\quad \quad  	For the third point the proof of the raising steps theorem
 gives that we have again  $\displaystyle v_{0}:=\sum_{j=1}^{N}{\chi
 _{j}u_{j}}$  with  $\displaystyle u_{j}\in \Lambda ^{\beta }(\Omega
 \cap U_{j})$  because the distance in  ${\mathbb{C}}^{n}$  and
 in  $\displaystyle \Omega \cap U_{j}$  are equivalent and the
  $\varphi _{j}$  are biholomorphic, as in remark~\ref{raisingSteps212},
 hence the  $\displaystyle \varphi _{j}$   send H\"older classes
 to the same H\"older classes.\ \par 
\quad \quad  	Then we apply the H\"older part of Kerzman's result~\cite{Kerzman71} :\ \par 
\quad \quad   $\displaystyle \forall \varphi \in L_{(0,1)}^{\infty }(\Omega
 ),\ \bar \partial \varphi =0,\ \exists w\in L^{\infty }(\Omega
 )::\bar \partial w=\omega \ and\ \forall \alpha <1/2,\ w\in
 \Lambda ^{\alpha }(\Omega ),\ {\left\Vert{w}\right\Vert}_{\Lambda
 ^{\alpha }(\Omega )}\leq C_{\alpha }{\left\Vert{\varphi }\right\Vert}_{L^{\infty
 }(\Omega )}\ ;$ \ \par 
 i.e. we choose  $\displaystyle r_{0}=\infty ,$  and we have
 a correction  $\displaystyle w\in \Lambda ^{\alpha }$  for any
  $\displaystyle \alpha <1/2.$  Hence again   $\bar \partial
 (v_{0}-w)=\omega $  and  $\displaystyle v_{0}-w\in \Lambda ^{\beta
 }$  because  $\displaystyle \beta <1/2$  and we choose of course
  $\displaystyle \alpha =\beta .$   $\blacksquare $ \ \par 

\subsection{Case of  $\displaystyle (p,q)$  forms by use of
 Demailly-Laurent's estimates.}
\quad \quad  	We can also remove the fact that  $\omega $  must have compact
 support when  $\displaystyle r\geq 2$  by using a theorem of
 Demailly-Laurent(~\cite{DemaillyLaurent87}, Remarque 4, page
 596) but the price to paid here is that the manifold has to
 be equipped with a metric with null curvature, in order to avoid
 parasitic terms.\ \par 
\quad \quad  	The proof is identical to the previous one and we get\ \par 
\begin{Theorem}
 Let  $\displaystyle \Omega $  be a  $\displaystyle s.p.c.$ 
 domain in the Stein manifold  $\displaystyle X$  equipped with
 a metric with null curvature. \par 
\quad  $\bullet $ 	 If  $\displaystyle 1\leq r<2(n+1)$  and  $\omega
 \in L^{r}_{(p,q)}(\Omega ),\ \bar \partial \omega =0.$  Then
 there is a  $\displaystyle (p,q-1)$  current  $u$  such that\par 
\quad \quad \quad \quad \quad 	 $\bar \partial u=\omega $  and  $\displaystyle u\in L^{s}_{(p,q-1)}(\Omega
 ),\ {\left\Vert{u}\right\Vert}_{s}\leq C{\left\Vert{\omega }\right\Vert}_{L_{(p,q)}^{r}(\Omega
 )}$  with  $\displaystyle \ \frac{1}{s}=\frac{1}{r}-\frac{1}{2(n+1)}.$ \par 
\quad \quad 	 $\bullet $  If  $\displaystyle r=2(n+1)$   and  $\omega \in
 L^{r}_{(p,q)}(\Omega ),\ \bar \partial \omega =0.$  Then for
 any  $\displaystyle s<\infty $  there is a  $\displaystyle (p,q-1)$
  current  $u$  such that\par 
\quad \quad \quad \quad \quad 	 $\bar \partial u=\omega $  and  $\displaystyle u\in L^{s}_{(p,q-1)}(\Omega
 ),\ {\left\Vert{u}\right\Vert}_{s}\leq C{\left\Vert{\omega }\right\Vert}_{L_{(p,q)}^{r}(\Omega
 )}.$ \par 
\end{Theorem}

\section{On convex domains of finite type.}
\quad \quad  	As another application of this method we get the following theorem :\ \par 
\begin{Theorem}
 Let  $\displaystyle \Omega $  be a domain in  ${\mathbb{C}}^{n}$
  such that locally around any point of  $\displaystyle \partial
 \Omega ,$   $\displaystyle \Omega $  is biholomorphic to a convex
 domain of finite type at most of type  $m\ ;$  then if  $\omega
 \in L^{r}_{(p,q)}(\Omega ),\ \bar \partial \omega =0$  if  $\displaystyle
 1\leq r\leq 2,$  there is a  $\displaystyle (p,q-1)$  current
  $u$  such that, with  $\displaystyle \gamma :=\min (\frac{1}{mn+2},\frac{1}{r}-\frac{1}{2}),$
  and  $\displaystyle \ \frac{1}{s}=\frac{1}{r}-\gamma ,$ \par 
\quad \quad \quad \quad \quad 	 $\bar \partial u=\omega $  and  $\displaystyle u\in L^{s}_{(p,q-1)}(\Omega
 ),\ {\left\Vert{u}\right\Vert}_{L_{(p,q-1)}^{s}(\Omega )}\leq
 C{\left\Vert{\omega }\right\Vert}_{L_{(p,q)}^{r}(\Omega )}.$ \par 
\end{Theorem}
\quad \quad  	Proof.\ \par 
We know that near any point  $\zeta $  of the boundary of  $\displaystyle
 \Omega $  we have that  $\displaystyle \Omega $  is biholomorphic
 to a  convex domain of finite type  $\displaystyle D_{\zeta }.$ \ \par 
\quad \quad  	In the case of convex domains of finite type, K. Diederich,
 B. Fischer and J-E. Fornaess~\cite{DieFisFor99}, A. Cumenge~\cite{Cumeng01}
 and B. Fischer~\cite{Fischer01} proved a theorem completely
 analogous to theorem~\ref{rS13}, with the gain given by  $\displaystyle
 \ \frac{1}{s}=\frac{1}{r}-\frac{1}{mn+2},$  so applying it we
 can solve  $\bar \partial $  with  $\displaystyle L^{r}-L^{s}$
  estimates with  $\displaystyle \ \frac{1}{s}=\frac{1}{r}-\frac{1}{mn+2}$
  for  $\displaystyle n<mn+2$  in  $D_{\zeta }$  and the biholomorphism
 gives the same result near  $\zeta $  in  $\displaystyle \Omega .$ \ \par 
Hence if  $\omega $  is a  $\displaystyle (p,q)$  form in  $\displaystyle
 L^{r}(\Omega ),$  there is a local solution  $u$  of  $\bar
 \partial u=\omega $  in  $\displaystyle L^{s}(D\cap \Omega )$
  where  $D$  is biholomorphic to a convex domain of finite type.
 On the other hand, the domain  $\displaystyle \Omega $  is pseudo
 convex, because locally biholomorphic to a pseudo convex one,
 and bounded, so we can use H\"ormander's theorem which says
 that we can solve globally the  $\bar \partial $  equation with
  $\displaystyle L^{2}$  estimates in such a domain  $\displaystyle
 \Omega .$ \ \par 
Hence we can apply the raising steps theorem~\cite{AmarSt13}
 to get a global  $\displaystyle L^{r}-L^{s}$  solution in  $\displaystyle
 \Omega .$   $\blacksquare $ \ \par 
\quad \quad  	To deal with  $\displaystyle r\geq 2$  we can apply exactly
 the same method as for theorem~\ref{lrEstimees33} to get\ \par 
\begin{Theorem}
Let  $\displaystyle \Omega $  be a domain in  ${\mathbb{C}}^{n}$
  such that locally around any point of  $\displaystyle \partial
 \Omega ,$   $\displaystyle \Omega $  is biholomorphic to a convex
 domain of finite type at most of type  $m\ ;$  Let  $\displaystyle
 \ 2<r\leq mn+2$  and  $\omega \in L^{r,c}_{(p,q)}(\Omega ),\
 \bar \partial \omega =0$  if  $\displaystyle 1\leq q<n$  and
  $\displaystyle \omega \in L^{r}_{(p,q)}(\Omega ),\ \omega \perp
 {\mathcal{H}}_{p}^{r'}(\Omega )$  if  $\displaystyle q=n\ ;$
  then there is a  $\displaystyle (p,q-1)$  current  $u$  such
 that, with  $\displaystyle \ \frac{1}{s}=\frac{1}{r}-\frac{1}{mn+2},$ \par 
\quad \quad \quad \quad \quad 	 $\bar \partial u=\omega $  and  $\displaystyle u\in L^{s}_{(p,q-1)}(\Omega
 ),\ {\left\Vert{u}\right\Vert}_{s}\leq C{\left\Vert{\omega }\right\Vert}_{r}.$
 \par 
\end{Theorem}
\ \par 
\ \par 
\ \par 

\bibliographystyle{C:/texlive/2012/texmf-dist/bibtex/bst/base/plain}

\end{document}